**Computing all $S$-integral points on elliptic curves**

Attila Pethő [†], Horst G. Zimmer, Josef Gebel and Emanuel Herrmann[‡]

# 1 Introduction

Let $a, b \in \mathbb{Z}$ with $\Delta_0 = 4a^3 + 27b^2 \neq 0$. Define a finite set of places $S = \{q_1, \ldots, q_{s-1}, q_s = \infty\}$ of $\mathbb{Q}$ and put $Q = \max\{q_1, \ldots, q_{s-1}\}$. Consider the group $E(\mathbb{Q})$ of rational points or Mordell-Weil group of the elliptic curve

$$E: \quad y^2 = x^3 + ax + b \tag{1}$$

over $\mathbb{Q}$ with discriminant $\Delta_0$ and absolute invariant $j = 12^3 \frac{4a^3}{\Delta_0}$. The multiplicative height of a rational point $P = (x, y) \in E(\mathbb{Q})$ is defined as the following product over all places $q$ of $\mathbb{Q}$ (including $q = \infty$):

$$H(P) = \prod_q \max\{1, |x|_q\},$$

where the $|x|_q$'s are the normalized multiplicative absolute values of $\mathbb{Q}$ corresponding to the places $q$ and satisfying the product formula.
Let $E(\mathbb{Z}_S)$ denote the set of $S-$integral points of $E(\mathbb{Q})$, i.e.

$$E(\mathbb{Z}_S) = \{P = (x, y) \in E(\mathbb{Q}) | H_S(P) \leq 1\},$$

where

$$H_S(P) = \prod_{q \notin S} \max\{1, |x|_q\}.$$

Siegel [22] proved in 1929 that the number of integral points on an elliptic curve $E$ over an algebraic number field $\mathbb{K}$ is finite, and Mahler [18] generalized this result in 1934 to $S$-integral points. Moreover, beginning with the pioneering work of Baker [1], several authors derived bounds for the size of the coordinates of integer points on elliptic curves $E$ over $\mathbb{K}$. Coates [5] proved that $E(\mathbb{Z}_S)$ contains only finitely many elements and that these can be effectively computed. The proofs used algebraic number-theoretical methods; the equation (1) was transformed into finitely many $S-$unit equations.
A completely different method for proving the finiteness of $E(\mathbb{Z}_S)$ was proposed by S. Lang [16] and Zagier [31]. Here one makes use of the group structure of $E(\mathbb{Q})$ and invokes properties of ordinary and $q$-adic elliptic logarithms. After S. David [7] established an explicit lower bound for linear forms in complex elliptic logarithms, Lang's idea could be transformed into an algorithm for computing all integer points on elliptic curves. This was done by Gebel, Pethő and Zimmer [10] and independently by Stroeker and Tzanakis [27]. The algorithm of Gebel,

[†]Research partially supported by Hungarian National Foundation for Scientific Research Grant No. 16791 and 16975.
[‡]Research partially supported by the Siemens AG.

Pethő and Zimmer was implemented by Gebel in the computer algebra package SIMATH and is available free of charge.
N. Smart [25] generalized the above method to an algorithm for finding all S-integral points on elliptic curves. Unfortunately, his approach, which seems to work well in practice, depends on an unproved lower bound for linear forms in $q$-adic elliptic logarithms. If the rank of $E(\mathbb{Q})$ is at most two, then such a bound was obtained by Rémond and Urfels [21] and applied by Gebel, Pethő and Zimmer [11] to find all S-integral points on Mordell's curves $y^2 = x^3 + k$, with $|k| \leq 10^4$ and such that the rank of the curve is strictly less then three.
In this note we combine the advantages of the methods of Siegel-Baker-Coates and of Lang-Zagier to overcome the absence of an explicit lower bound for linear forms in $q$-adic elliptic logarithms. This was made possible by a recent, completely explicit upper bound for the $S-$integral solutions of (1) established by Hajdu and Herendi [15].
In Section 2 we state our theoretical result and in the subsequent section we prove it. In Section 4 we deal with the computational aspects of the main result and in Section 5 we explain how we compute $q$-adic elliptic logarithms. Finally we show through an example how our method works in practice.

## 2 Notations and the Theorem

To state our result we require several definitions. Let $P_1, \ldots, P_r$ denote a basis of the Mordell-Weil group and let $g$ be the order of the torsion subgroup $E(\mathbb{Q})_{tors}$ of $E(\mathbb{Q})$. We have $g \leq 12$ by a famous result of Mazur [20]. Let $\hat{h}$ denote the Néron-Tate height as well as the associated symmetric bilinear form on $E(\mathbb{Q})$. Designate by $\lambda$ the smallest eigenvalue of the positive definite regulator matrix $(\hat{h}(P_i, P_j))_{1 \leq i,j \leq r}$.
Let $\wp(u)$ be the Weierstraß $\wp-$function corresponding to the curve $E(\mathbb{C})$. Let $\Omega$ be its fundamental lattice and $\omega$ its real period. There exists, for any $P = (x, y) \in E(\mathbb{C})$, an element $u \in \mathbb{C}/\Omega$ such that $(x, y) = (\wp(u), \wp'(u))$. This $u$ is called the *(complex) elliptic logarithm* of $P$. In the sequel $u_{i,\infty}$ denotes the elliptic logarithm of $P_i, i = 1, \ldots, r$. We put $u'_{i,\infty} = gu_{i,\infty}/\omega$.
For a prime $q$ of $\mathbb{Q}$, $\tilde{E}$ denotes the reduced curve $E$ modulo $q$. Let $\mathcal{N}_q = \#\tilde{E}(\mathbb{F}_q)$ be the number of rational points on $\tilde{E}/\mathbb{F}_q$. With the order $g$ of the torsion group, we define the number

$$m = m_q = \mathrm{lcm}(g, \mathcal{N}_q).$$

Finally, for the finite places $q \in S$, let $u'_{i,q}$ denote the $q-$adic elliptic logarithm of $mP_i$ for $i = 1, \ldots, r$. Now we are in a position to state our main result.

**Theorem** *Assume that the $S-$integral point $P = (x, y) \in E(\mathbb{Z}_S)$ has the representation*

$$P = \sum_{i=1}^{r} n_i P_i + T \tag{2}$$

*with* $n_i \in \mathbb{Z}, i = 1, \ldots, r$, *and* $T$ *a torsion point of* $E(\mathbb{Q})$. *For* $N = \max\{|n_i|, i = 1, \ldots, r\}$, *we have*

$$N \leq N_0 = \sqrt{(c_1/2 + c_2)/\lambda}, \tag{3}$$

*with*

$$c_2 = \log \max\{(8|a|)^{1/4}, (32|b|)^{1/6}\},$$

$$c_1 = 7 \cdot 10^{38s+49} s^{20s+15} Q^{24} (\log^* Q)^{4s-2} c_3 (\log c_3)^2 (c_3 + 20(s-1)c_3 + \log(ec_4)),$$

*where* $\log^* Q = \max\{\log Q, 1\}$ *for* $Q = \max\{q_1, \ldots, q_{s-1}\}$,

$$c_3 = \frac{32}{3}\sqrt{|\Delta_0|}\left(8 + \frac{1}{2}\log|\Delta_0|\right)^4,$$

$$c_4 = 10^4 \max\{16a^2, 256\sqrt{|\Delta_0|}^3\}$$

*with* $\Delta_0 = 4a^3 + 27b^2$.
*Moreover, there exists a place* $q \in S$ *such that*

$$\left|\sum_{i=1}^{r} n_i u'_{i,q} + n_{r+1}\right|_q \leq c_5 \exp\{-(\lambda/s)N^2 + c_2/s\}, \tag{4}$$

*with* $n_{r+1} \in \mathbb{Z}$ *if* $q = \infty$ *and* $n_{r+1} = 0$ *otherwise, and with* $c_5 = \sqrt{8}g/\omega$ *if* $q = \infty$ *and* $c_5 = 1$ *otherwise.*

**Remarks**

1. If $r \leq 2$, then we can use the estimate of Rémond and Urfels [21] for linear forms in $q$-adic elliptic logarithms to derive an upper bound for $N$ (see [11]). The example shown in Section 6 reveals that this bound can be much smaller than our bound $N_0$.

2. The proof of the Theorem will show that, in the Proposition below, $c_1$ can be replaced by any upper bound for the height of $S-$integer points in $E(\mathbb{Q})$. For example, in the case of $s = 1$, one can take the smaller constant

$$c'_1 = 5 \cdot 10^{64} c_3 \log c_3 (c_3 + \log c_4))$$

   instead of specializing to $s = 1$ in the formula for $c_1$ (see [15], Theorem 1). Then we obtain

$$N \leq N'_0 = \sqrt{(c'_1/2 + c_2)/\lambda}$$

   instead of (3).

In the case of $s = 1$, Gebel, Pethő and Zimmer [10] * as well as Stroeker and Tzanakis [27] established another bound for $N$, namely

$$N \leq N_1 = 2^{r+3}\sqrt{c_6 c_7} \log^{(r+3)/2}(c_7(r+3)^{r+3})$$

*The constants given here differ from those in [10] because in [10] (on pages 179-180) two corrections are necessary: First $r$ must be replaced by $r+1$ and second, $C$ requires an adjustment to the new version of the Theorem in [7].

with

$$
\begin{aligned}
c_6 &= \max\left\{\frac{\log(2\sqrt{2}\sqrt[3]{4}g/\omega)}{\lambda}, 1\right\}, \\
c_7 &= \max\left\{\frac{C}{\lambda}, 10^9\right\}\left(\frac{h}{2}\right)^{r+2}\prod_{i=0}^{r}\log V_i, \\
C &= 2.9\cdot 10^{6(r+2)}4^{2(r+1)^2}(r+2)^{2r^2+13r+23.3}, \\
h &= \log\max\{4|aj_2|, 4|bj_2|, |j_1|\} \text{ for } j=\frac{j_1}{j_2}, \\
& \quad j_1, j_2 \in \mathbb{Z}, \gcd(j_1, j_2) = 1,
\end{aligned}
$$

and

$$
\begin{aligned}
&\log V_i \geq \max\left\{\hat{h}(P_i), h, \tfrac{3\pi|u'_{i,\infty}|_\infty^2}{Im\tau}\right\}, \quad i = 1, \ldots, r, \\
&\log V_0 \geq \max\{h, \tfrac{3\pi}{\text{Im }\tau}\}.
\end{aligned}
$$

It seems hard to compare the size of $N_0'$ and $N_1$ because they depend on different parameters of the curve. The numerical examples of Section 6 suggest that $N_1 \leq N_0'$ if $r \leq 2$, but $N_1 > N_0'$ if $r > 2$. Therefore, it makes sense to compute both bounds $N_0'$ and $N_1$ before starting the reduction procedure, which is based on inequality (4) and does not depend on the way, by which an upper bound for $N$ was obtained.

# 3 Proof of the Theorem

The proof of the Theorem is divided into three parts: First we obtain an explicit upper bound for the ordinary height $h(P)$ of the $S-$integral points $P$ in $E(\mathbb{Q})$. In the second part we prove a relation between $h(P)$ and $N$ and derive (3). Finally, we prove (4) by using properties of complex and $q-$adic elliptic logarithms.

## 3.1 Estimation of $h(P)$.

**Proposition** *Let $E(\mathbb{Q})$ be the elliptic curve (1) and $S$ be the finite set of rational places (containing the infinite place $\infty$) defined above. Then all $S-$integral points*

$$
P = (x, y) = \left(\frac{\xi}{\zeta^2}, \frac{\eta}{\zeta^3}\right), \quad \xi, \eta, \zeta \in \mathbb{Z}, \quad \gcd(\xi, \zeta) = 1 = \gcd(\eta, \zeta)
$$

*in $E(\mathbb{Q})$ satisfy*

$$
\max\{h(x), h(y)\} \leq c_1,
$$

*with the constant $c_1$ defined in the Theorem, where*

$$
h\left(\frac{\alpha}{\beta}\right) = \log\max\{|\alpha|, |\beta|\}, \quad \alpha, \beta \in \mathbb{Z}, \gcd(\alpha, \beta) = 1.
$$

**Proof**: [15], Theorem 2.

### 3.2 Height estimates and the proof of (3).

Let us start out from an elliptic curve

$$E: \; Y^2 + a_1XY + a_3Y = X^3 + a_2X^2 + a_4X + a_6 \qquad (a_i \in \mathbb{Z})$$

in long Weierstraß form over $\mathbb{Q}$. Tate's quantities are

$$\begin{aligned} &b_2 = a_1^2 + 4a_2,\; b_4 = 2a_4 + a_1a_3,\; b_6 = a_3^2 + 4a_6, \\ &b_8 = a_1^2a_6 + 4a_2a_6 - a_1a_3a_4 + a_2a_3^2 - a_4^2, \\ &c_4 = b_2^2 - 24b_4,\; c_6 = -b_2^3 + 36b_2b_4 - 216b_6. \end{aligned}$$

$E$ has discriminant

$$\Delta = -b_2^2b_8 - 8b_4^3 - 27b_6^2 + 9b_2b_4b_6 = -16\Delta_0$$

and absolute invariant

$$j = \frac{c_4^3}{\Delta} = \frac{j_1}{j_2} \qquad (j_1, j_2 \in \mathbb{Z},\, \gcd(j_1, j_2) = 1).$$

The invariant differential is

$$\omega = \frac{dx}{2y + a_1x + a_3} = \frac{dy}{3x^2 + 2a_2x + a_4 - a_1y}.$$

We begin by estimating heights. In the introduction we defined the multiplicative height of a rational point $P = (x, y) \in E(\mathbb{Q})$. The *ordinary additive height* is then

$$h(P) := \frac{1}{2}\,\log\, H(P)$$

and the *Néron-Tate* or *canonical additive height* is

$$\hat{h}(P) := \lim_{n\to\infty} \frac{h(2^nP)}{2^{2n}}.$$

From Theorem 5 in [34] (cf. also [24], [32], [33]), one readily derives the inequality

$$h(P) \geq \hat{h}(P) - c_2' \quad \text{ for } P \in E(\mathbb{Q})$$

with the constant

$$c_2' := \frac{1}{2}(\log\, 2 + \mu_\infty),$$

where

$$\mu_\infty := \log\, \max\{|b_2|_\infty, |b_4|_\infty^{\frac{1}{2}}, |b_6|_\infty^{\frac{1}{3}}, |b_8|_\infty^{\frac{1}{4}}\}.$$

We now specialize this height estimation to the curve in short Weierstrass form (1). We have $a_1 = a_2 = a_3 = 0$ and $a_4 = a, a_6 = b$, so that $\mu_\infty$ becomes

$$\mu_\infty = \log\max\{(2|a|)^{1/2}, (4|b|)^{1/3}, |a|^{1/2}\} = \log\max\{(2|a|)^{1/2}, (4|b|)^{1/3}\}.$$

Hence

$$h(P) \geq \hat{h}(P) - c_2 \quad \text{ for } P \in E(\mathbb{Q}). \tag{5}$$

with the constant $c_2$ defined in the Theorem.
Let, as before, $P_1, \dots, P_r$ be a basis of the free part of the Mordell-Weil group $E(\mathbb{Q})$. Then $P \in E(\mathbb{Q})$ has the (unique) representation (2) with $n_i \in \mathbb{Z}, i = 1, \dots r$, and $T$ a torsion point. Put $N = \max\{|n_i|, i = 1, \dots, r\}$. The canonical height $\hat{h}$ can be extended to the $r-$dimensional real vector space $E(\mathbb{Q}) \otimes_{\mathbb{Z}} \mathbb{R}$ and the extension of $\hat{h}$ is a positive definite quadratic form on this space.
As explained in [10], this leads to the lower estimate on non-torsion points

$$\hat{h}(P) \geq \lambda N^2$$

involving the smallest eigenvalue $\lambda$ of the regulator matrix associated with the basis points $P_1, \dots, P_r \in E(\mathbb{Q})$ via the canonical height $\hat{h}$. This, together with (5), implies that

$$h(P) \geq \lambda N^2 - c_2 \quad \text{for } P \in E(\mathbb{Q}) \setminus E(\mathbb{Q})_{tors}. \tag{6}$$

Let now $P = (x, y) \in E(\mathbb{Z}_S)$. Then, since $h(P) = \frac{1}{2} \log h(x)$, we obtain

$$h(P) \leq c_1/2$$

by the Proposition, with $c_1$ as defined in the Theorem. On combining the lower and upper estimates for $h(P)$, we obtain

$$\lambda N^2 - c_2 \leq c_1/2.$$

This implies the asserted inequality (3):

$$N \leq \sqrt{(c_1/2 + c_2)/\lambda)}$$

### 3.3 Elliptic logarithms and the proof of inequality (4).

Let now $P = (x, y) \in E(\mathbb{Q})$ be an $S$-integral point and choose $q \in S$ such that

$$|x|_q = \max\{|x|_{q_1}, \dots, |x|_{q_{s-1}}, |x|_\infty\}. \tag{7}$$

In the sequel we keep this $q$ fixed. First we conclude that

$$H(P) \leq |x|_q^s \qquad \text{for } s := \sharp S$$

and hence that

$$h(P) \leq \frac{s}{2} \log |x|_q. \tag{8}$$

Combining (6) and (8) yields (cf. [8])

$$\frac{1}{|x|_q^{1/2}} \leq c_8 \exp{(-c_9 N^2)} \tag{9}$$

with

$$c_8 := \exp\left(\frac{c_2}{s}\right), \ c_9 := \frac{\lambda}{s}.$$

To transform the upper bound for $|x|_q^{-\frac{1}{2}}$ into the upper bound for the left hand side of the inequality (4), we must study the properties of elliptic logarithms. Two cases are to be distinguished, the complex case of $q = \infty$ and the $q$-adic case.

Case 1: $q = \infty \in S$.

Let $\wp(u)$ be the Weierstraß function associated to the curve (1). Let $\Omega = \langle \omega, \omega_2 \rangle$ be the period lattice of $\wp(u)$, generated by the two fundamental periods $\omega$ and $\omega_2$, where $\omega > 0$ is real and $\omega_2$ is complex. Let $P = (x, y) \in E(\mathbb{Q})$ and denote by $u = u(P)$ its elliptic logarithm, i.e. that residue class $u \in \mathbb{C} \bmod \Omega$ for which $P = (\wp(u), \wp'(u))$.

Then we have

$$u \equiv \int_x^{\infty} \frac{dt}{\sqrt{t^3 + at + b}} \pmod{\Omega}$$

(see [31]).

Let

$$t^3 + at + b = (t - \alpha)(t - \beta)(t - \gamma)$$

be the decomposition of the cubic polynomial into linear factors with roots

$$\alpha < \beta < \gamma \quad \text{if } \alpha, \beta, \gamma \in \mathbb{R}$$

and

$$\alpha = \bar{\beta} \in \mathbb{C}, \quad \gamma \in \mathbb{R} \text{ otherwise.}$$

Define

$$M = \begin{cases} 0 & \text{, if } \quad \gamma \geq 0 \\ \dfrac{\exp \mu_\infty}{2^{1/3} - 1}, & \text{if } \quad \gamma < 0 \end{cases}$$

and put

$$x_0 = \max\{\alpha + \beta, 2\gamma\} + M.$$

With this notation, we have the

**Lemma** *Let $x \in \mathbb{R}$ such that $x > \max\{0, x_0\}$. Then*

$$\int_x^{\infty} \frac{dt}{\sqrt{t^3 + at + b}} \leq \frac{\sqrt{8}}{\sqrt{x}}.$$

**Proof:** [10], pp. 176 and 181-183. Note that in [10], the number $x$ is assumed to be the first coordinate of an integer point of $E(\mathbb{Q})$, but in the proof this assumption was not needed.

Let now $P = (x, y)$ be an $S-$integral point on $E(\mathbb{Q})$ with first coordinate of absolute value

$$|x|_\infty > \max\{x_0, \exp \mu_\infty\}. \tag{10}$$

Denote by $u$ its elliptic logarithm. Assume furthermore that (7) holds with $q = \infty$. Then, by the Lemma and by (9), we obtain, on normalizing $u$ to $0 < |u|_\infty \leq 1/2\omega$,

$$|u|_\infty \leq \frac{\sqrt{8}}{|x|_\infty^{1/2}} \leq \sqrt{8}\, c_8\, \exp(-c_9 N^2). \tag{11}$$

Let $u_{i,\infty} = u(P_i), i = 1, \ldots, r$, denote the elliptic logarithms of $P_1, \ldots, P_r$. Then

$$u(gP) = u\left(\sum_{i=1}^{r} n_i(gP_i)\right) = \sum_{i=1}^{r} n_i(gu_{i,\infty}),$$

where $g$ is the order of the torsion subgroup of $E(\mathbb{Q})$. Moreover, there exists an integer $n_{r+1} \in \mathbb{Z}$ such that

$$|n_1(gu_{1,\infty}) + \ldots + n_r(gu_{r,\infty}) + n_{r+1}\omega|_\infty \leq \sqrt{8}\ g\ c_8\ \exp(-c_9 N^2).$$

Dividing by $\omega$ we obtain (4) in the case of $q = \infty$.

Case 2: $q = q_i \in S$ (for some $i \in \mathbb{N}$ such that $1 \leq i \leq s-1$).
Here we use $q$-adic elliptic logarithms (cf. [8], [23], [25]).

For the convenience of the reader, we explain in some detail how one proceeds in the $q$-adic case. Let $\mathbb{Q}_q$ be the $q$-adic completion of $\mathbb{Q}$ and $\mathbb{Z}_q$ be its ring of $q$-adic integers. Denote by

$$E_1(\mathbb{Q}_q) := \{P \in E(\mathbb{Q}_q)\ \mid\ \tilde{P} = \tilde{\mathcal{O}}\}$$

the kernel of the reduction map modulo $q$, where $E$ is regarded as a curve over $\mathbb{Q}_q$ and $\tilde{P}$, $\tilde{\mathcal{O}}$ are the reduced points $P, \mathcal{O}$ modulo $q$. Designate by $\mathcal{E}(q\mathbb{Z}_q)$ the formal group associated to $E$ (cf. [23], [25]). We consider the isomorphism

$$\begin{array}{ccc} \mathcal{E}(q\mathbb{Z}_q) & \longrightarrow & E_1(\mathbb{Q}_q) \\ z & \longmapsto & \left\{ \begin{array}{cc} 0 & \text{if } z = 0 \\ (\frac{z}{w(z)}, -\frac{1}{w(z)}) & \text{if } z \neq 0 \end{array} \right\}, \end{array}$$

where

$$z = -\frac{x}{y},\ w(z) = -\frac{1}{y}.$$

The equation for $w = w(z)$ inferred from the long Weierstraß equation for $E(\mathbb{Q})$ becomes

$$w = z^3 + (a_1 z + a_2 z^2)w + (a_3 + a_4 z)w^2 + a_6 w^3 = f(z, w).$$

A recursive procedure based on this equation (see [23]) leads to the power series

$$\begin{aligned} w = \ & z^3 + a_1 z^4 + (a_1^2 + a_2)z^5 + (a_1^3 + 2a_1 a_2 + a_3)z^6 \\ & + (a_1^4 + 3a_1^2 a_2 + 3a_1 a_3 + a_2^2 + a_4)z^7 + \cdots \\ & \in \mathbb{Z}[a_1, a_2, a_3, a_4, a_6][[z]]. \end{aligned}$$

This is the unique power series in $z$ satisfying the relation

$$w(z) = f(z, w(z)).$$

From it we also get the Laurent series for $x$ and $y$, viz.

$$\begin{aligned} x(z) &= \tfrac{z}{w(z)} &= \tfrac{1}{z^2} - \tfrac{a_1}{z} - a_2 - a_3 z - (a_4 + a_1 a_3)z^2 - \cdots, \\ y(z) &= -\tfrac{1}{w(z)} &= -\tfrac{1}{z^3} + \tfrac{a_1}{z^2} + \tfrac{a_2}{z} + a_3 + (a_4 + a_1 a_3)z + \cdots. \end{aligned} \tag{12}$$

The invariant differential has the expansion

$$\begin{aligned}\omega(z) = \quad & (1 + a_1 z + (a_1^2 + a_2)z^2 + (a_1^3 + 2a_1a_2 + a_3)z^3 \\ & +(a_1^4 + 3a_1^2a_2 + 6a_1a_3 + a_2^2 + 2a_4)z^4 + \cdots)dz.\end{aligned}$$

Note that in these expansions the coefficients of the powers of $z$ each have the same weight depending on the exponent of $z$.

The $q$-adic elliptic logarithm is now the image under the homomorphism to the additive group $\hat{G}_a$ (over the completion $\mathbb{C}_q$ of the algebraic closure of $\mathbb{Q}_q$.) defined as follows:

$$\begin{array}{lll}\psi_q : E_1(\mathbb{Q}_q) & \longrightarrow & \hat{G}_a \\ P = (x, y) & \longmapsto & \psi_q(P) = \int \omega(z) = z + \frac{d_2}{2}z^2 + \frac{d_3}{3}z^3 + \cdots .\end{array}$$

In particular, the $q$-adic logarithm $\psi_q$ has the properties

$$\psi_q(P + Q) = \psi_q(P) + \psi_q(Q)$$

and

$$|\psi_q(P)|_q = |z|_q = \left|-\frac{x}{y}\right|_q .$$

As before, let $\tilde{E}$ be the reduced curve $E$ modulo $q$, and denote by $\mathcal{N}_q = \sharp\tilde{E}(\mathbb{F}_q)$ the number of rational points on $\tilde{E}/\mathbb{F}_q$. With the number

$$m = \mathrm{lcm}(g, \mathcal{N}_q)$$

defined in Section 2, we have from the Lutz filtration of $E$ (see [17]):

$$mP_i =: P_i' \in E_1(\mathbb{Q}_q) \quad (i = 1, \ldots, r)$$

for the generating points $P_i$ of $E(\mathbb{Q})$, and

$$mP_{r+1} = \mathcal{O}$$

for the torsion points $P_{r+1} \in E(\mathbb{Q})_{tors}$.

The representation (2) of an $S$-integral point $P = (x, y) \in E(\mathbb{Q})$ then gives rise to the representation

$$P' = n_1'P_1 + \cdots + n_r'P_r = n_1P_1' + \cdots n_rP_r' \quad (n_i' := mn_i \in \mathbb{Z}) \tag{13}$$

of its $m$-multiple $P' = (x', y') = mP \in E_1(\mathbb{Q}_q)$. In analogy to (12), we have the Laurent series

$$x' = \frac{z'}{w(z')} = \frac{1}{z'^2} - \frac{a_1}{z'} - a_2 - a_3z' - (a_4 + a_1a_3)z'^2 - \cdots ,$$

and this expansion entails the estimate

$$|x'|_q \leq \frac{1}{|z'|_q^2} = \frac{1}{|u_q'|_q^2}, \tag{14}$$

where we use the abbreviating notation $u'_q := \psi_q(P')$ for the elliptic logarithm of $P'$.

On combining the inequalities (9) and (14) and observing that $|x'|_q \geq |x|_q$, we end up with the upper estimate for the $q$-adic elliptic logarithm $u'_q = \psi_q(P')$ of the point $P' = (x', y') = mP = m(x, y)$:

$$|u'_q|_q \leq \frac{1}{|x'|_q^{1/2}} \leq \frac{1}{|x|_q^{1/2}} \leq c_8 \ \exp(-c_9 N^2).$$

From (13), we have the relation

$$u'_q = n'_1 u_{1,q} + \cdots + n'_r u_{r,q} = n_1 u'_{1,q} + \ldots + n_r u'_{r,q} \quad (n'_i = m n_i \in \mathbb{Z})$$

between the elliptic logarithms $u'_q = \psi_q(P')$ of the $m$-multiple $P' = mP$ of the given $S$-integral point $P$, $u_{i,q} = \psi_q(P_i)$ of the generating points $P_i$ and $u'_{i,q} = \psi_q(P'_i)$ of their $m$-multiples $P'_i = mP_i \in E(\mathbb{Q})$, and inequality (4) is proved in this case too.

# 4 Computational aspects of the Theorem.

The Theorem, together with some computational techniques for elliptic curves and with diophantine approximation methods, in practice makes it possible to compute all $S-$integral solutions of (1). To this end, we must carry out the following steps:

**(i)** Compute a basis $P_1, \ldots, P_r$ of the Mordell-Weil group $E(\mathbb{Q})$. This can be achieved, for example, by employing the algorithm [9] (see also [19], [33]). Having obtained the basis, one can easily compute the constants $c_1, c_2$ and $\lambda$ and thus the upper bound $N_0$ for $N$, defined by the inequality (3).

**(ii)** $N_0$ is usually very large (see the examples), but we also have the estimate (4) at one of the places $q \in S$. Some observations on the computation of rational approximations of the elliptic logarithms with the necessary very high precision will be made in the next section. The solutions of the system of inequalities (3) and (4) are then found by the reduction method of de Weger [29]. Since this method is well-described in [25], we omit the details here. Of course, one must perform de Weger's reduction for all places $q \in S$. Let $\mathcal{M}_q$ denote the bound for $N$ after the reduction with respect to a prime $q \in S$, assuming that (7) holds for this $q \in S$. Then we change the value of $N_0$ to $\max\{\mathcal{M}_q; q \in S\}$. If possible we iterate the reduction.

**(iii)** Experience shows that the last value of $N_0$ obtained after several reductions is sufficiently small in order to enable one to test the points $\sum_{i=1}^r n_i P_i + T$ within the range $N = \max\{|n_i|, i = 1, \ldots, r\} \leq N_0$ and with $T \in E_{tors}(\mathbb{Q})$ for $S$-integrality.

However, this can be a lengthy procedure. Therefore, we describe here some ideas as to how to proceed in a more economical way. If the $r-$tuple $\underline{n} = (n_1, \ldots, n_r)$ of integer vectors is such that $N = \max\{|n_i|, i = 1, \ldots, r\} \leq N_0$, then we first test for $q \in S$ by floating point approximation (if $q = \infty$) or by $q-$adic approximation (if $q$ is a finite place) of $u'_{i,q}$, whether or not (4) is true. If, for an $\underline{n}$, inequality (4) does not hold for some $q \in S$, then there is no $S-$integral point on the curve given by the vector $\underline{n}$ in the representation (2). Otherwise we must test all linear combinations $\sum_{i=1}^{r} n_i P_i + T$ with varying torsion points $T$ for $S-$integrality.

**(iv)** There is one extra search to be carried out eventually. This is due to the fact that, if (7) holds for $q = \infty$, then (11) and (4) hold only for those $S-$integral points $P = (x, y)$, whose first coordinate satisfies the inequality (10). We must also determine the $S$-integral points for which (10) is violated. Then, by (7), we have

$$|x|_{q_i} \leq \max\{x_0, \exp \mu_\infty\}, \quad i = 1, \ldots, s-1.$$

Write

$$x = \frac{\xi}{q_1^{\alpha_1} \cdots q_{s-1}^{\alpha_{s-1}}} \quad \text{with } \gcd(\xi, q_i) = 1, \quad i = 1, \ldots, s-1. \tag{15}$$

This yields

$$|x|_{q_i} = q_i^{\alpha_i},$$

so that

$$\alpha_i \leq \log \max\{x_0, \exp \mu_\infty\} / \log q_i, \quad i = 1, \ldots, s-1. \tag{16}$$

Hence the exponents $\alpha_i$ are bounded. We note that the exponents $\alpha_i$ are even. The negation of inequality (10), together with the boundedness relations (16) for the coordinates (15), shows that there are only finitely many $S$-integral points, which are to be tested in this case.

# 5 Computation of the $q$-adic elliptic logarithms

The complex elliptic logarithms can be computed either by applying the arithmetic-geometric mean method (see H.Cohen [6]) or by employing a method of Zagier [31].

In contrast to the case of complex elliptic logarithms, there exists, for the computation of $q$-adic elliptic logarithms, a quadratically convergent procedure only for Tate curves (see [14]) *. That is why we must apply a different $q$-adic procedure here. In fact, we simply use the formulae of Section 3. However, we must overcome some technical difficulties if we need to calculate the $q$-adic elliptic logarithms to a precision of up to 1000 or more digits. Of course, we can compute the formal Laurent series expansion of $w(T)$ by Newton iteration. However, unfortunately the coefficients in this expansion soon become very large.

*We wish to thank J.H.Silverman for calling our attention to the paper [14]

We circumvent this problem by invoking the Lutz filtration of $E$ over $\mathbb{Q}_q$ and calculating $q$-adic approximations to a sufficient high precision as follows. Suppose that, for some given positive integer $t$, we have computed an approximation

$$\tilde{\psi}(T) := \sum_{i=1}^{t} \frac{d_i}{i} T^i \equiv \psi(T) \pmod{T^{t+1}}$$

of the formal elliptic integral $\psi(T) = \int w(T)\, dT$, where $d_1 = 1$. We want to compute an approximation $\widehat{\psi_q(P)}$ of the $q$-adic elliptic logarithm $\psi_q(P)$ of $P \in E_1(\mathbb{Q}_q)$ such that

$$|\psi_q(P) - \widehat{\psi_q(P)}|_q \leq q^{-n}$$

for some positive integer $n$. (In practice, $t$ is fixed and $n$ depends on $q$ and is much larger than $t$.) Put $V = [\frac{n}{t}]$ and define $P_V = q^V P$. Then we have $P_V \in E_{V+1}(\mathbb{Q}_q)$ by the Lutz filtration theorem [17]. Now we can compute a suitable $q$-adic approximation $\tilde{P}_V$ of the point $P_V$ since $E$ is considered as a curve over $\mathbb{Q}_q$ as in Section 3.3. Here a suitable approximation means that, for the point $P_V = (x_{P_V}, y_{P_V}) \in E(\mathbb{Q}_q)$, we choose $\tilde{P}_V = (\tilde{x}_{P_V}/q^{2(V+1)}, \tilde{y}_{P_V}/q^{3(V+1)})$ with $\tilde{x}_{P_V}, \tilde{y}_{P_V} \in \mathbb{Z}$ such that

$$\left| \frac{x_{P_V}}{y_{P_V}} - \frac{\tilde{x}_{P_V}}{\tilde{y}_{P_V}} q^{V+1} \right|_q \leq q^{-n'},$$

with $n' = n + V$. To this end we take $\tilde{z}_{P_V} \in \mathbb{Z}$ such that

$$|\tilde{z}_{P_V}|_q = \left| -\frac{\tilde{x}_{P_V}}{\tilde{y}_{P_V}} q^{V+1} \right|_q \quad \text{and} \quad \left| \tilde{z}_{P_V} + \frac{\tilde{x}_{P_V}}{\tilde{y}_{P_V}} q^{V+1} \right|_q \leq q^{-n'}.$$

Then we have

$$\left| \tilde{z}_{P_V} + \frac{x_{P_V}}{y_{P_V}} \right|_q \leq q^{-n'},$$

which implies

$$\begin{aligned} \left| \tilde{\psi}(\tilde{z}_{P_V}) - \tilde{\psi}\left( -\frac{x_{P_V}}{y_{P_V}} \right) \right|_q &= \left| \sum_{i=1}^{t} \frac{d_i}{i} \left( \tilde{z}_{P_V}^i - \left( -\frac{x_{P_V}}{y_{P_V}} \right)^i \right) \right|_q \\ &\leq \left| \tilde{z}_{P_V} + \frac{x_{P_V}}{y_{P_V}} \right|_q \leq q^{-n'}. \end{aligned}$$

On the other hand, we have

$$\begin{aligned} &\left| \psi\left( -\frac{x_{P_V}}{y_{P_V}} \right) - \tilde{\psi}\left( -\frac{x_{P_V}}{y_{P_V}} \right) \right|_q = \left| \sum_{i=t+1}^{\infty} \frac{d_i}{i} \left( -\frac{x_{P_V}}{y_{P_V}} \right)^i \right|_q \\ &\leq \left| \frac{d_{t+1}}{t+1} \left( -\frac{x_{P_V}}{y_{P_V}} \right)^{t+1} \right|_q \leq \frac{q^{-(t+1)(V+1)}}{|t+1|_q} \leq q^{-n'} \end{aligned}$$

by the choice of $V$. Thus

$$\left|\psi\left(-\frac{x_{P_V}}{y_{P_V}}\right)-\tilde{\psi}(\tilde{z}_{P_V})\right|_q \le q^{-n'}.$$

We have

$$\psi\left(-\frac{x_{P_V}}{y_{P_V}}\right)=\psi_q(P_V)=q^V\psi_q(P),$$

hence

$$|\psi_q(P)-\tilde{\psi}(\tilde{z}_{P_V})/q^V|_q \le q^{-n'+V}=q^{-n},$$

i.e. $\tilde{\psi}(\tilde{z}_{P_V})/q^V$ is the required approximation of $\psi_q(P)$, the $q$-adic elliptic logarithm of $P$.

# 6 Comparison of the Theorem with earlier estimates

We first consider the $S=\{3,5,\infty\}-$integer points of the curve

$$E_{Sm}(\mathbb{Q}) : y^2+y=x^3+x^2-2x,$$

treated by Smart [25]. He proved that $P_1=(0,0), P_2=(1,0)$ form a basis of the Mordell-Weil group of $E_{Sm}$ over $\mathbb{Q}$ and that $E_{Sm}$ has no non-trivial torsion points. Moreover, he showed that if

$$P=n_1P_1+n_2P_2$$

is an $S-$integral point on $E_{Sm}$ with $\max\{|n_1|,|n_2|\}\le 10^{40}$, then $\pm(n_1,n_2)=$ $(0,1),(0,2),(0,3),(1,0),(1,1),(1,-1),(1,-4),(2,0),(2,-1),(2,2),(2,-2),$ $(2,-3),(3,0),(3,1),(3,-2),(4,-1),(1,2)$.
Using the transformation $X=36x+12, Y=216y+108$, we obtain the short Weierstraß form

$$E'_{Sm}(\mathbb{Q}) : Y^2=X^3-3024X+46224.$$

It is clear that if $P=(x,y)$ is an $S-$integral point on $E_{Sm}$, then $P'=(x',y')=(36x+12,216y+108)$ is an $S-$integral point on $E'_{Sm}$. Thus it suffices to find all $S-$integral points on $E'_{Sm}$. Let $P'=(x',y')$ be an $S-$integral point on $E'_{Sm}$ with

$$P'=n_1P'_1+n_2P'_2$$

for $P'_1=(12,108), P'_2=(48,108)$. Then the Theorem implies that

$$N=\max\{|n_1|,|n_2|\}\le N_0=6\cdot 10^{121},$$

as follows by a simple computation. This bound is much larger then the bound used by Smart. However, by performing one de Weger reduction with this *unconditional* bound, we obtain

$$\max\{|n_1|,|n_2|\}\le 106.$$

This shows that Smart indeed found all $S-$integral points on $E_{Sm}$.

As the rank of $E_{Sm}$ is two, one can also apply the theorem of Rémond and Urfels [21] to estimate $N$ in the way described in [11]. Using that method one easily verifies that

$$N \leq 4.7 \cdot 10^{41}.$$

Here the worst case appears when $\max\{|u|_3, |u|_5, |u|_\infty\} = |u|_\infty$. This is a much smaller bound for $N$ than the one which follows from the Theorem, and hence we obtain another proof that Smart found all $S-$integral points on $E_{Sm}$, this time without any further computation.

In Remark (2) and thereafter, we gave two different bounds $N_0'$ and $N_1$ for $N$ in the case $s = 1$. Now we are going to compare them for some curves, for which all integer points were found by the elliptic logarithm method.

The first set of examples is taken from [13], where the classical Mordell equations

$$y^2 = x^3 + k$$

were solved † for $|k| \leq 10^5$. In Table 1 we list, for some $k$'s, the rank $r_k$ of the curve $y^2 = x^3 + k$ and the values of $N_1$ and $N_0'$.

**Table 1**

| $k$ | $r_k$ | $N_1$ | $N_0'$ |
|---|---|---|---|
| 108 | 1 | $2.8 \cdot 10^{26}$ | $1.8 \cdot 10^{41}$ |
| 225 | 2 | $1.3 \cdot 10^{41}$ | $4.5 \cdot 10^{41}$ |
| 1025 | 3 | $5.5 \cdot 10^{60}$ | $3.1 \cdot 10^{42}$ |
| 2089 | 4 | $1.1 \cdot 10^{84}$ | $7.7 \cdot 10^{42}$ |
| -28279 | 5 | $2.1 \cdot 10^{112}$ | $2.0 \cdot 10^{42}$ |

Next, in Table 2, we take some examples from [3], where the family of elliptic curves

$$y^2 = x^3 - 36x - 864k(k-1)(2k-1)$$

was considered.

**Table 2**

| $k$ | $r_k$ | $N_1$ | $N_0'$ |
|---|---|---|---|
| 1 | 1 | $8.9 \cdot 10^{23}$ | $1.3 \cdot 10^{41}$ |
| 3 | 2 | $5.8 \cdot 10^{39}$ | $1.8 \cdot 10^{44}$ |
| 7 | 3 | $2.4 \cdot 10^{60}$ | $5.9 \cdot 10^{45}$ |
| 20 | 4 | $2.1 \cdot 10^{86}$ | $2.9 \cdot 10^{47}$ |

†The rank one cases left open in [13] have been settled by K. Wildanger in his PhD thesis "Über das Lösen von Einheiten- und Indexformgleichungen in algebraischen Zahlkörpern mit einer Anwendung auf die Bestimmung aller ganzen Punkte einer Mordellschen Kurve", Berlin 1997. In those cases, the Mordell curves have no integral points - as conjectured by the authors.

Finally, we consider four curves from [8], [10] and from [28] of the form

$$y^2 = x^3 + ax + b.$$

In Table 3 we exhibit, for certain coefficients $a, b$, the rank $r$ of the curve defined by (1) and the corresponding values of $N_1$ and $N_0'$.

**Table 3**

| $a$ | $b$ | $r$ | $N_1$ | $N_0'$ |
|---|---|---|---|---|
| -203472 | 18487440 | 5 | $2.3 \cdot 10^{111}$ | $7.9 \cdot 10^{47}$ |
| -1642032 | 628747920 | 6 | $1.1 \cdot 10^{144}$ | $2.1 \cdot 10^{49}$ |
| -147952 | 21183760 | 7 | $2.7 \cdot 10^{187}$ | $1.1 \cdot 10^{47}$ |
| -5818216808130 | 5401285759982786436 | 8 | $8.67 \cdot 10^{224}$ | $2.33 \cdot 10^{58}$ |

The examples indicate, that $N_1 < N_0'$, if $r \le 2$, but $N_1 > N_0'$, if $r > 2$. Of course, for rank 0 curves it is quite easy to find all integer points by examining the algebra of the curve. Our comparison shows that $N_1$ is significantly smaller than $N_0'$ if $r = 1$. For $r = 2$, $N_1$ and $N_0'$ are of a similar size, but $N_1$ appears to be slightly smaller than $N_0'$. For $r \ge 3$, the constant $N_0'$ becomes much smaller than $N_1$, and with growing rank $r$, the difference becomes rather significant. This can be very well observed in the above rank 7 and 8 examples. The reason for this phenomenon probably lies in the fact that in David's [7] lower bound, which enters into the bound $N_1$, the constant

$$C = 2.9 \cdot 10^{6(r+2)} 4^{2(r+1)^2} (r+2)^{2r^2+13r+23.3}$$

shows up, and this constant depends exponentially on $r^2$. On the other hand, the constant term of $N_0'$ does not depend on $r$.
Our examples make it plausible that, in order to find all integer points on elliptic curves, it makes sense to first compute both bounds $N_1$ and $N_0'$ and then start the reduction procedure with the smaller one.

# 7 An example

To illustrate our method, we choosen the rank-four curve

$$E: \quad y^2 = x^3 - 172x + 505 \tag{17}$$

from the paper A. Wiman [30] (see also Bremner and Tzanakis [4]). Wiman showed that there are at least 58 integer points on $E$. Here we prove that there are exactly 58 integer points and, moreover, 144 $S$- integer points on $E$, where $S = \{3, 5, 7, \infty\}$. More precisely, we prove the following theorem.
**Theorem′** *Let $S = \{3, 5, 7, \infty\}$. Then all $S$-integer points $P = (x, \pm y)$ with $x = \frac{\xi}{\zeta^2}, y = \frac{\eta}{\zeta^3}; \xi, \eta, \zeta \in \mathbb{Z}, \gcd(\xi, \zeta) = \gcd(\eta, \zeta) = 1$ on $E$ over $\mathbb{Q}$ are those displayed in Table 5.*
**Proof.** First we compute the quantities $c_1, c_2, c_3$ and $c_4$ of the Theorem in Section 2 to obtain $N_0$ in dependence on $\lambda$ only. For this purpose, we do not need any geometrical information about the elliptic curve $E$. We have $a = -172$ and $b = 505$, hence

$$\Delta_0 = 2198992, \quad c_2 = 1.81, \quad c_3 = 8.7 \cdot 10^8 \quad \text{and} \quad c_4 = 8.35 \cdot 10^{15}.$$

Taking $s = 4$ and $Q = 7$, we obtain

$$c_1 = 4.564 \cdot 10^{305},$$

and this yields

$$N_0 = 4.777 \cdot 10^{152}/\sqrt{\lambda}.$$

Now we consider $E$ as an elliptic curve over $\mathbb{Q}$. Inserting the coefficients of $E$ into SIMATH we obtain, after some minutes of computation, that $E$ has trivial torsion, its rank is 4, and a basis of its Mordell-Weyl group is given by the points:

$$P_1 = (12, 13), P_2 = (-14, 13), P_3 = (-1, 26) \quad \text{and} \quad P_4 = (38, 221).$$

The regulator of the curve is $R = 2.79532$, the real period is $\omega = 0.808974$ and the least eigenvalue corresponding to the chosen basis is $\lambda = 0.7467531$.
We have $g = 1$, hence $m_q = \mathcal{N}_q$, and the actual values of $m_q$ are $m_3 = 7, m_5 = 10$ and $m_7 = 12$. By the Theorem, for each $q \in S$, we have to solve the following diophantine approximation problems

$$\begin{aligned} \left| \sum_{i=1}^{4} n_i u'_{i,q} + n_{r+1} \right|_q &\leq \exp\{-0.149 \cdot N^2 + c_{10}\}, \\ N = \max\{|n_1|, \ldots, |n_4|\} \leq N_0 &= 5.53 \cdot 10^{152} \end{aligned}$$

for each $q \in S$, where $c_{10} = 1.62$, if $q = \infty$ and $c_{10} = 0.362$ otherwise.
In order to reduce the huge upper bound for $N$, we first take $q = \infty$ and perform de Weger reduction with $C = 10^{910}$. We obtain the new upper bound $N \leq \mathcal{M}_\infty = 69$. Next, for each $q \in S \setminus \{\infty\}$, we compute the $q$-adic elliptic logarithms of $P_i, i = 1, \ldots, 4$, with precision at least

$$n_3 = 1281, \;\; n_5 = 875, \;\; n_7 = 723.$$

This precision is necessary in order to carry out the $q$-adic de Weger reduction. For this purpose, we use the method of Section 5. We have chosen $t = 200$ and, for safety reasons,

$$V_3 = 8, \;\; V_5 = 6, \;\; V_7 = 5,$$

and computed the approximations of $q^{V_q} P_i$ for $q \in \{3, 5, 7\}$ and $i = 1, \ldots, 4$ with

$$n'_3 = 1700, \;\; n'_5 = 1300, \;\; n'_7 = 1100$$

$q$-adic digit numbers. In Table 4 we have listed the first and last three significant digits of $\psi_q(m_q P_i)$, i.e., the digits $a_2, a_3, a_4, a_{n_q-2}, a_{n_q-1}$ and $a_{n_q}$ assuming, that

$$\psi_q(m_q P_i) \approx \sum_{i=1}^{n_q} a_i q^i.$$

Note that, as $m_q P_i \in E_1(\mathbb{Q})$, we always have $a_1 = 0$.

**Table 4**

| $i$ | 1 | 2 | 3 | 4 |
|---|---|---|---|---|
| $P_i$ | $(12, 13)$ | $(-14, 13)$ | $(-1, 26)$ | $(38, 221)$ |
| $\psi_3(7P_i)$ | $(1, 2, 1, \dots, 0, 2, 2)$ | $(2, 0, 1, \dots, 2, 1, 2)$ | $(1, 1, 0, \dots, 2, 2, 2)$ | $(0, 2, 2, \dots, 0, 1, 1)$ |
| $\psi_5(10P_i)$ | $(3, 0, 4, \dots, 3, 3, 4)$ | $(4, 4, 0, \dots, 0, 3, 1)$ | $(2, 2, 3, \dots, 3, 3, 0)$ | $(4, 3, 4, \dots, 0, 1, 1)$ |
| $\psi_7(12P_i)$ | $(0, 6, 0, \dots, 2, 3, 4)$ | $(5, 2, 5, \dots, 0, 2, 4)$ | $(2, 1, 1, \dots, 5, 6, 2)$ | $(1, 4, 4, \dots, 1, 2, 1)$ |

Now we perform the $q$-adic de Weger reduction with the values $C_3 = 3^{1281}, C_5 = 5^{875}$ and $C_7 = 7^{723}$ and obtain the new bound

$$N_0 = 124 = \max\{\mathcal{M}_\infty, \mathcal{M}_3, \mathcal{M}_5, \mathcal{M}_7\}.$$

This new upper bound for $N$ can be further reduced. On repeating this reduction process 3-times, we eventually get $N \leq 17$, which cannot be reduced any further. In order to compute all $S$-integral solutions of (17) with this last bound, we use the sieving procedure described in [8]. Finally, we must also determine the $S$-integral points for which (10) is violated. In the case under consideration, we have

$$|x|_{q_i} < \max\{22.56, 31.86\} = 31.86, \quad i = 1, \dots, 4, \tag{18}$$

thus $\alpha_2 \leq 4$ and $\alpha_3, \alpha_5 \leq 2$, where we used the notation of Section 4, part (iv). A simple computation shows that all the solutions of (17) satisfying the above inequality (18) are among those already found. Theorem′ is proved. □

**Table 5**

$S$-integral points $P = (x, y) = \left(\frac{\xi}{\zeta^2}, \frac{\eta}{\zeta^3}\right) = \sum_{i=1}^{4} n_i P_i$ on

$E: \; y^2 = x^3 - 172x + 505$

for $S = \{3, 5, 7, \infty\}$

| | |
|---|---|
| rank | 4 |
| basis | $P_1 = (12, 13), \; P_2 = (-14, 13), \; P_3 = (-1, 26), \; P_4 = (38, 221)$ |
| torsion | trivial |

| # | $\xi$ | $\eta$ | $\zeta$ | $F$ | $(n_1, n_2, n_3, n_4)$ |
|---|---|---|---|---|---|
| 1 | $-14$ | 13 | 1 | | $(0, 1, 0, 0)$ |
| 2 | $-12$ | 29 | 1 | | $(0, 0, 1, 1)$ |
| 3 | $-10$ | 35 | 1 | | $(-1, 0, -1, 0)$ |
| 4 | $-8$ | 37 | 1 | | $(0, -1, 0, -1)$ |
| 5 | $-2$ | 29 | 1 | | $(1, 0, -1, 0)$ |

| $E: \; y^2 = x^3 - 172x + 505$ | (continued) |
|---|---|

| # | $\xi$ | $\eta$ | $\zeta$ | $F$ | $(n_1, n_2, n_3, n_4)$ |
|---|---|---|---|---|---|
| 6 | $-1$ | 26 | 1 | | $(0, 0, 1, 0)$ |
| 7 | 2 | 13 | 1 | | $(-1, -1, 0, 0)$ |
| 8 | 3 | 4 | 1 | | $(-1, 1, 0, 0)$ |
| 9 | 12 | 13 | 1 | | $(1, 0, 0, 0)$ |
| 10 | 14 | 29 | 1 | | $(-1, 0, 0, -1)$ |
| 11 | 16 | 43 | 1 | | $(0, -1, -1, 0)$ |
| 12 | 19 | 64 | 1 | | $(1, 1, 1, 1)$ |
| 13 | 24 | 101 | 1 | | $(0, 1, -1, 0)$ |
| 14 | 31 | 158 | 1 | | $(1, 0, 0, -1)$ |
| 15 | 38 | 221 | 1 | | $(0, 0, 0, 1)$ |
| 16 | 76 | 653 | 1 | | $(-2, 0, 0, 0)$ |
| 17 | 90 | 845 | 1 | | $(0, -1, -1, -1)$ |
| 18 | 131 | 1492 | 1 | | $(-1, -1, 1, -1)$ |
| 19 | 168 | 2171 | 1 | | $(1, 1, 1, 0)$ |
| 20 | 284 | 4781 | 1 | | $(0, -2, 0, 0)$ |
| 21 | 415 | 8450 | 1 | | $(2, 0, 0, 1)$ |
| 22 | 467 | 10088 | 1 | | $(0, 1, -1, -1)$ |
| 23 | 1046 | 33827 | 1 | | $(1, 0, 2, 1)$ |
| 24 | 1266 | 45043 | 1 | | $(-1, 1, 1, -1)$ |
| 25 | 1314 | 47629 | 1 | | $(-1, 1, -1, 1)$ |
| 26 | 3028 | 166621 | 1 | | $(1, 0, -2, 0)$ |
| 27 | 9502 | 926237 | 1 | | $(-2, -2, 0, -1)$ |
| 28 | 11283 | 1198496 | 1 | | $(-2, 1, 1, 1)$ |
| 29 | 1402464 | 1660877429 | 1 | | $(0, -3, -1, -2)$ |
| 30 | $-128$ | 233 | 3 | 3 | $(-2, -1, 0, 1-)$ |
| 31 | 4 | 559 | 3 | 3 | $(-1, 0, -1, -1)$ |
| 32 | 28 | 1 | 3 | 3 | $(0, 0, 1, -1)$ |
| 33 | 106 | 287 | 3 | 3 | $(-1, 1, 1, 1)$ |
| 34 | 160 | 1495 | 3 | 3 | $(-1, -1, 1, 0)$ |
| 35 | 823 | 23374 | 3 | 3 | $(0, 2, 0, 1)$ |
| 36 | 11134 | 1174769 | 3 | 3 | $(-1, -2, -2, -1)$ |
| 37 | 33524044 | 194104052639 | 3 | 3 | $(0, 0, 2, -2)$ |
| 38 | $-1085$ | 14680 | 9 | $3^2$ | $(-1, 1, 2, 0)$ |
| 39 | $-926$ | 22789 | 9 | $3^2$ | $(2, -1, 0, 0)$ |
| 40 | $-536$ | 26819 | 9 | $3^2$ | $(-1, -1, 0, 1)$ |
| 41 | 4612 | 305227 | 9 | $3^2$ | $(-1, 0, -2, 0)$ |
| 42 | 2182 | 811433 | 27 | $3^3$ | $(-1, -2, -1, -2)$ |
| 43 | 14119 | 13113463 | 27 | $3^3$ | $(-1, 2, 0, -1)$ |
| 44 | 1709021 | 62200333 | 81 | $3^4$ | $(2, 0, 2, 1)$ |
| 45 | $-181$ | 4628 | 5 | 5 | $(2, 1, 0, 0)$ |

| $E:\ y^2 = x^3 - 172x + 505$ | (continued) |
|---|---|

| # | $\xi$ | $\eta$ | $\zeta$ | $F$ | $(n_1, n_2, n_3, n_4)$ |
|---|---|---|---|---|---|
| 46 | 7824 | 691457 | 5 | 5 | $(-2, -1, -1, -2)$ |
| 47 | 1064 | 33137 | 5 | 5 | $(2, -1, -1, 0)$ |
| 48 | 294 | 1303 | 5 | 5 | $(-1, -2, 0, -1)$ |
| 49 | 3314 | 189863 | 5 | 5 | $(-1, -1, -1, 1)$ |
| 50 | −324 | 2951 | 5 | 5 | $(1, 0, -1, 1)$ |
| 51 | 1191 | 39614 | 5 | 5 | $(-1, 0, 2, 1)$ |
| 52 | 346 | 3481 | 5 | 5 | $(1, -1, 1, 1)$ |
| 53 | −64 | 3809 | 5 | 5 | $(-1, 1, 0, 1)$ |
| 54 | −7034 | 497861 | 25 | $5^2$ | $(0, 0, -1, -2)$ |
| 55 | 1754 | 104117 | 25 | $5^2$ | $(0, 1, -2, 0)$ |
| 56 | 7405599 | 15552621382 | 625 | $5^4$ | $(1, -2, -2, 1)$ |
| 57 | 334 | 53677 | 15 | $3 \times 5$ | $(0, 2, 1, 0)$ |
| 58 | 441046 | 292380031 | 45 | $3^2 \times 5$ | $(3, 0, 0, -1)$ |
| 59 | −26 | 12571 | 7 | 7 | $(-2, 0, 1, 0)$ |
| 60 | −657 | 6866 | 7 | 7 | $(-1, -2, -1, 0)$ |
| 61 | 16 | 7267 | 7 | 7 | $(1, -1, 0, 1)$ |
| 62 | 128 | 2941 | 7 | 7 | $(0, 1, 2, 1)$ |
| 63 | 554 | 811 | 7 | 7 | $(1, 1, -1, -1)$ |
| 64 | 562 | 2197 | 7 | 7 | $(0, 1, -1, 1)$ |
| 65 | 1082 | 29653 | 7 | 7 | $(2, -1, -1, -1)$ |
| 66 | 140260 | 52528645 | 7 | 7 | $(1, 0, 0, -2)$ |
| 67 | −33666 | 1488019 | 49 | $7^2$ | $(0, -2, -1, -2)$ |
| 68 | −24097 | 4109846 | 49 | $7^2$ | $(-1, 0, 1, -2)$ |
| 69 | 399988 | 252199571 | 49 | $7^2$ | $(1, -2, -2, 0)$ |
| 70 | −2498 | 333593 | 21 | $3 \times 7$ | $(-1, 2, -1, 0)$ |
| 71 | 28596 | 4114769 | 35 | $5 \times 7$ | $(2, 0, 2, 0)$ |
| 72 | 430259 | 157313848 | 175 | $5^2 \times 7$ | $(-1, 3, 1, 1)$ |

Attila Pethő
Institut of Mathematics and Informatics
Kossuth Lajos Universität
H-4010 Debrecen, P.O.Box 12
Hungary

Horst G. Zimmer
Josef Gebel
Emanuel Herrmann
Fachbereich 9 Mathematik
Universität des Saarlandes
Postfach 151150
D-66041 Saarbrücken
Germany